\documentclass{MR3}

\usepackage{lipsum}

\usepackage{url,hyperref,bold-extra}
\usepackage{amssymb,amsfonts,amsmath}
\usepackage{pifont,yfonts,mathrsfs}
\usepackage{latexsym,bm}

\newenvironment{prf}
{\noindent\mbox{\textbf{\textsc{Proof}}:}}{\hfill{\mbox{\ding{113}}}\bigskip}

\begin{document}
\title{A Non-Constructive Proof of Cantor's Theorem}
\titlerunning{A Non-Constructive Proof of Cantor's Theorem} %
\author{\sc \textbf{Saeed~Salehi}} 
% ...
\authorrunning{S.\ Salehi} 
\affiliation{%
 Plaksha University, IT City Rd, Sector 101A, \\ Mohali, Punjab 140306, India. \textsf{root@saeedsalehi.ir}
   }

\MakeTitle

Cantor showed that there are {\em hierarchies of infinities}, not just {\em one infinity}; this gave a new life to the concept of infinity by bringing it into the domain of mathematics and freeing it from philosophers' dominance.

\bigskip

\noindent
{\textbf{\textsc{Theorem}}} ({\sc Georg Cantor}, 1891).

\noindent
{\em There can be no surjection from a set $A$ onto its powerset $\mathscr{P}(A)$.}

\bigskip

Cantor's original proof uses his celebrated
diagonal argument, showing  that for a  function $f \colon A \rightarrow \mathscr{P}(A)$, the anti-diagonal set $$\mathcal{D}_0 = \{a \in A\mid
a \not\in f(a)\}$$ is not in the range of $f$. Some alternative proofs are constructive, like Cantor's, as they explicitly construct some sets that are not in the range of $f$:
 If the binary relation $\mathcal{R} \subseteq A^2$ is defined by $x\mathcal{R}y
 {\text{ if and only if } } y \in f(x)$, for $x,y \in A$, and  $n > 0$ is an integer, the set $$\mathcal{D}_n = \{a \in A\mid \boldsymbol\not \exists \{x_i\}_{i=1}^{n}  :  a\mathcal{R}x_1\mathcal{R}\cdots\mathcal{R}x_n
\mathcal{R}a\}$$ is out of the $f$'s range   (cf. \S24 of {\bf \cite{Quine}}); so is the set $$\mathcal{D}_\infty = \{
a \in A\mid \boldsymbol\not \exists
\{x_i\}_{i=1}^{\infty}  : a\mathcal{R}x_1\mathcal{R}x_2\mathcal{R}\cdots\}$$ (see {\bf \cite{Raja}}).
There is a  non-constructive proof that
shows the nonexistence of an injection $h \colon \mathscr{P}(A) \rightarrow A$ (see {\bf \cite{Boolos}}); a part of the proof is constructive as it builds two subsets $B,C \subseteq A$ such that   $h(B) = h(C)$ but $B \neq C$ (one can even have   $B\subsetneq C$ and $h(B) = h(C) \in C \setminus B$, see \S3 of {\bf \cite{KS}}).
However, that this   implies the nonexistence of a surjection $f \colon A \rightarrow \mathscr{P}(A)$ requires the Axiom of Choice, thus making the whole argument non-constructive.
We give another non-constructive proof (which seems new as it is not listed in the ``various proofs'' of {\bf \cite{Raja2}}).

\bigskip

\begin{prf}

\noindent
For finite sets, this follows from the Pigeonhole principle: if $A$ has $n$ elements, then prove by induction (on $n$) that $\mathscr{P}(A)$ has $2^n$ elements and that $2^n > n$ holds.
If $A$ is not finite, then partition it into some finite subsets, such as $A = \bigcup_{i\in I}A_i$, where the  $A_i$s, for $i \in I$, are nonempty and pairwise disjoint.
We show that for a given function $f \colon A \rightarrow \mathscr{P}(A)$, there exists some $B \subseteq A$ that is out of $f$'s range without explicitly describing what it could look like.
For each $i \in I$, consider $f_i \colon A_i \rightarrow \mathscr{P}(A_i)$, defined by $f_i(x) = f(x) \cap A_i$, for $x \in A_i$. There exists a subset $B_i \subseteq A_i$ that is not in the range of $f_i$; recall that each $A_i$ is finite.  We show that the set $B = \bigcup_{i\in I}B_i$ is not in the range of $f$. If, otherwise, $B = f(\alpha)$ holds for some $\alpha \in A$, then there is a unique $\kappa \in I$ such that $\alpha \in A_\kappa$. Now, we have  $B_\kappa = B \cap A_\kappa = f(\alpha) \cap A_\kappa
 = f_\kappa(\alpha)$,  but this contradicts the choice of $B_\kappa$ (which was supposed to be out of the  $f_\kappa$'s range).
\end{prf}

Let us notice the use of the Axiom of Choice in the above proof, once in partitioning the set $A$ into $\{A_i\}_{i\in I}$ and once in {\em choosing} the subsets $\{B_i\}_{i\in I}$ of $A_i$s.
If one takes the $A_i$s to be singletons, then one gets Cantor's anti-diagonal set $\mathcal{D}_0 = \bigcup_{a\in A}[\{a\} \setminus f(a)]$, since the only subset of the singleton $\{a\}$ that is not in the range of $g \colon \{a\} \rightarrow \mathscr{P}(\{a\})$ is the set  $\{a\} \setminus g(a)$.
A non-constructive proof is obtained, when
the finite subsets $A_i$s of $A$ are taken to have more than one element.



\begin{references}{}







\bibitem{Boolos}
{\sc G. Boolos}, {\em Constructing Cantorian counterexamples}, \textbf{\textit{Journal of  Philosophical Logic}}   26:3 (1997),  237--239 (with an editorial note by  V.~McGee).
{\sc doi}: \href{https://doi.org/10.1023/A:1004209106100}{\tt 10.1023/A:1004209106100}.
Reprinted in:
\textbf{\textit{Logic, Logic, and Logic}},  R.~Jeffrey (ed.), Harvard University Press (1998), pp.~339--341.
{\sc jStor}: \href{https://www.jstor.org/stable/30227093}{\tt 30227093}

\bibitem{KS}
{\sc A. Karimi \& S. Salehi}, {\em Diagonal arguments and  fixed points}, \textbf{\textit{Bulletin of the  Iranian Mathematical  Society}} 43:5
 (2017),  1073--1088. \\
\url{http://bims.iranjournals.ir/article_979.html}



\bibitem{Quine}
{\sc W. Quine}, \textbf{\em Mathematical Logic}, Harvard University Press (revised 1951).
{\sc isbn}:~\href{https://isbnsearch.org/isbn/9780674554504}{\tt 9780674554504}



\bibitem{Raja}
{\sc N. Raja}, {\em A negation-free proof of Cantor's theorem}, \textbf{\textit{Notre Dame Journal of  Formal  Logic}} 46:2 (2005), 231--233.
{\sc doi}: \href{https://doi.org/10.1305/ndjfl/1117755152}{\tt 10.1305/ndjfl/1117755152}




\bibitem{Raja2}
{\sc N. Raja}, ``{\em Yet another proof of Cantor’s theorem}'', in: \textbf{\textit{Dimensions of Logical Concepts}}, J.-Y. B\'eziau, A. Costa-Leite (eds.),  Cole\c{c}\~{a}o CLE  54, Campinas (2009), pp.~209--217. \url{https://bit.ly/4lEpssR}




\end{references}
\end{document}